\font\bigtype=cmr12
\font\biggtype=cmbx12 scaled \magstep 2
\font\Smcps = cmcsc10 scaled \magstep 1

\long\def\title#1{\bgroup
\biggtype \parindent=0pt\rightskip = 0in plus 10in
\leftskip=\rightskip \spaceskip=.3333em \xspaceskip=.5em
\parfillskip=0pt \baselineskip=21pt #1

\egroup\vskip .4cm}

\newcount\sezione

\newcount\enunciato

\enunciato=1

\def\iniziosezione{\vskip0pt plus.3\vsize\penalty-50\vskip0pt
plus-.3\vsize\bigskip\bigskip\vskip\parskip}

\def\beginsectionnonumber#1\par{\message{#1}
\iniziosezione
\leftline{\Smcps#1}\nobreak\medskip}

\def\beginsection #1\par{\message{#1}\advance\sezione by
1\enunciato=1
\iniziosezione
\leftline{\Smcps\the\sezione\enspace#1}
\nobreak\medskip}

\outer\def\references{
\beginsection References

\parindent=0pt}

\def\w#1 {\def\quan{1.5cm}\smallbreak\noindent
\hangindent\quan\hangafter1\leavevmode\hbox to
\quan{[#1]\hfill}}

\def\proclaim #1. {\medbreak
\noindent{\bf#1.}
\enspace\sl}

\def\endproofsymbol{\clubsuit}

\def\proclaimr #1. {\proclaim #1. \rm}

\def\endproclaim{\advance\enunciato by
1\rm\par\ifdim\lastskip<\medskipamount
\removelastskip\penalty55\medskip\fi}

\outer\def\proof{\smallbreak\noindent{\bf Proof.}\enspace}

\outer\def\proofof #1. {\smallbreak\noindent{\bf Proof of #1.}\enspace}

\outer\def\endproof{~\hfill~$\endproofsymbol$\par
\ifdim\lastskip<
\medskipamount\removelastskip\medbreak\fi}

\def\H{{\rm H}}

\def\pic{\mathop{\rm Pic}}

\def\cl{\mathop{\rm Cl}}

\def\tr{\mathop{\rm tr}}

\def\into{\hookrightarrow}

\title{On the minimal normal compactification

of a polynomial in two variables}

\vskip.5cm

\centerline{{\bigtype Angelo Vistoli}\footnote *{Partially supported by the
University of Bologna, funds for selected research topics.}}\vskip .3cm
\centerline {Dipartimento di Matematica}
\centerline {Universit\`a di Bologna}
\centerline{Piazza di Porta San Donato 5}
\centerline{40127 Bologna, Italy}
\smallbreak
\centerline{E-mail address: \tt vistoli@dm.unibo.it}

\vskip 1cm

\beginsection Introduction

Let $C$ be an integral affine curve over a field $\kappa$,
$\alpha,\beta: C\into {\bf A}^2$ two closed embeddings. We say that $\alpha$
and $\beta$ are {\it equivalent\/} when there is an automorphism $\phi$ of
${\bf A}^2$ with $\phi\circ \alpha =
\beta$. It was stated by B. Segre ([Se]) and
proved by Suzuki in [Su] that when $\kappa = {\bf C}$ any embedding of
${\bf A}^1$ into ${\bf A}^2$ is equivalent to the standard embedding $t \mapsto
(t,0)$. This was generalized to the case that
$\kappa$ is arbitrary, and the degree of $f$ is prime to the characteristic
of $\kappa$, by Abhyankar and Moh ([AM]). This is
what is usually called the Abhyankar--Moh theorem. On the other hand, there
are many affine curves with an infinite number of non equivalent embeddings
into ${\bf A}^2$: for example, ${\bf A}^1\setminus\{0\}$.

Suzuki in [Su] also proves a very nice result: if $C$ is smooth and has only
one branch at infinity (that is, it is the complement
of a point in a smooth projective curve) and $f$ is a
generator of the ideal of $C$ in ${\bf A}^2$, then $C$ is an ordinary fiber
of $f$, that is,
$f$ is a topological fibration in a neighborhood of 0. In their important,
and arduous, article [AS] Abhyankar and Singh carry the study of this case
much further, over arbitrary fields; in particular, for example, such a curve
$C$ has at most finitely many nonequivalent embeddings, with appropriate
conditions on the characteristic of the field.

Another proof of Suzuki's theorem was given by Artal Bartolo, in [AB], based
on the results of [EN], relating knot theory with the theory of polynomials in
two variables. 

Now, let $f\colon {\bf A}^2\to
{\bf A}^1$ be a polynomial in two variables defined over an algebraically
closed field $\kappa$. We shall always assume that $f$ is primitive, that
is, that the generic fiber of $f$ is integral. We consider the {\it minimal
normal compactification} $\overline f
\colon X\to {\bf A}^1$ of
$f$, namely the only normal irreducible surface $X$ containing ${\bf A}^2$
as an open subset, together with a proper morphism $\overline f \colon X\to
{\bf A}^1$ extending $f$, with the property that each fiber of $f$ is dense
in the corresponding fiber of $\overline f$. It is often singular.

Let $E_1,\ldots,E_r$ be the horizontal components of $X\setminus{\bf A}^2$,
namely the irreducible components of $X\setminus{\bf A}^2$ that dominate
${\bf A}^1$. By standard results,  $E_1,\ldots,E_r$ are isomorphic to ${\bf
A}^1$ and do not intersect (Proposition~1). To each $E_i$ we associate two
integers. The first is the degree $e_i$ of $E_i$ over ${\bf A}^1$; one
can think of $e_1,\ldots,e_r$ as the orders of the orbits of the monodromy
group acting on the branches at infinity of a general fiber of $f$. The
second is the least positive integer $\delta_i$ such that $\delta_iE_i$ is
a Cartier divisor on $X$; since $E_i$ is smooth, we have that $\delta_i = 1$
if and only if $X$ has no singularities along $E_i$.

Our result (Theorem~1) says that if the characteristic of $\kappa$ is 0, the
greatest common divisor of
$\delta_1e_1,\ldots,\delta_re_r$ is 1. In particular, if there only one
component $E_1$, this maps isomorphically onto ${\bf A}^1$, and
$\overline f$ is smooth along $E_1$. So, if one of the fibers of
$f$ has only one branch at infinity, then there is a simultaneous resolution
of singularities at infinity of $f$. This easily implies the
Suzuki--Abhyankar--Moh embedding theorem.

One can show that the integers $e_i \delta_i$ coincide with the integers $m_i$
defined by Eisenbud and Neumann (see [AB], p.~102). So in characteristic~0 our
result follows from [EN], section~4, although our proof is shorter.

In characteristic $p$
we only get that the greatest common divisor of
$\delta_1e_1,\ldots,\delta_re_r$ is 1 when the degree of the polynomial is
prime to $p$ (see Theorem~2). This implies the Suzuki--Abhyankar--Moh embedding
theorem over a perfect field.

The proof of the Theorem~1 is entirely
straightforward, and very short; it uses standard topological methods, plus
some elementary facts about rational surface singularities. If one substitutes
ordinary topological cohomology with \'etale cohomology with ${\bf Z}_\ell $
coefficients, where $\ell$ is a prime different from the characteristic of
$\kappa$, one gets a proof of Theorem~2. We do not
include the proof of this general case, but anyone who is familiar with
\'etale cohomology will be able to reconstruct the details.

\beginsection Acknowledgments

I am grateful to Pierrette Cassou-Nogu\`es for several helpful comments about
the history of the results I quote. In particular she pointed out reference
[AB] to me.

\beginsection The results

Consider a complex polynomial in two variables, i.e., a morphism $f\colon
{\bf A}^2\to {\bf A}^1$ defined over {\bf C}. We shall always assume that
$f$ is primitive, that is, that $f$ is not constant, and not obtained by
composition with a polynomial in one variable of degree greater than 1. This
the same as saying that the generic fiber of $f$ is integral, or that the
subfield ${\bf C}(f)$ is algebraically closed in ${\bf C}(x,y)$.

We consider the minimal normal compactification $\overline f\colon X\to {\bf
A}^1$ of $f$, obtained by taking the closure $\Gamma$ of the graph of $f$ in
${\bf P}^2\times {\bf P}^1$, considering its normalization $X'$, and then
calling $X$ the inverse image of ${\bf A}^1$ in $X'$. Then $X$ is a normal
integral complex quasiprojective scheme over
$\kappa$, containing ${\bf A}^2$ as an open subscheme. Furthermore the
morphism $f$ extends to a morphism $\overline f \colon X\to {\bf A}^1$,
which has the useful property that every fiber of $f$ is dense inside the
corresponding fiber of $\overline f $. Let us call $E_1,\ldots,E_r$ the
irreducible components of the complement $E$ of ${\bf A}^2$ in $X$. Each of
the $E_1,\ldots,E_r$ is an affine integral curve dominating ${\bf A}^1$: we
will call $e_1,\ldots,e_r$ the degrees of $E_1,\ldots,E_r$ over ${\bf A}^1$.

Furthermore, the divisor class groups of the local rings of $X$ are finite,
because $X$ has rational singularities ([Li], Proposition 17.1.) We will call
$\delta_i$ the least common multiple of the orders of
$E_i$ in each of the divisor class groups of the local rings of $X$ at
points of $E_i$; clearly $\delta_i E_i$ is a Cartier divisor on $X$, while
$\delta E_i$ is not a Cartier divisor for any integer $\delta$ with $0 <
\delta < \delta_i$.

\proclaim Theorem 1.  Each of the $E_1,\ldots,E_r$ is isomorphic to
${\bf A}^1$, and they are pairwise disjoint.

Furthermore the greatest common divisor of the products
$\delta_1e_1,\ldots, \delta_re_r$ is 1.
\endproclaim

The first statement in the theorem is quite standard. It has an important
consequence; if $\delta_i$ is 1, that is, if $E_i$ is a Cartier divisor on
$X$, then $X$ is smooth at all point of $E_i$.

\proclaim Corollary 1.  Assume that $X$ has only one component at infinity.
Then all the fibers of $\overline f$ are integral and smooth at
infinity.

In particular, this happens when one of the fibers of $f$ has
only one branch at infinity.\endproclaim

This follows immediately from the theorem, because the hypothesis implies
that $\delta_1 = 1$, i.e., $X$ is smooth, and $e_1 = 1$, i.e., $E_1$
maps isomorphically onto ${\bf A}^1$.

{}From the corollary we get a new proof of the renowned
Suzuki--Abhyankar--Moh theorem. For this we only need to assume that
$\kappa$ is perfect.

\proclaim The Suzuki--Abhyankar--Moh theorem over {\bf C}.  Any embedding of
${\bf A}^1$ into ${\bf A}^2$ defined over {\bf C} is equivalent to the
standard embedding $t\mapsto(t,0)$.
\endproclaim

\proof.  Let $C$ be a curve in ${\bf
A}^2$ isomorphic to ${\bf A}^1$, $f\in {\bf C}[x,y]$ a generator of the
ideal of $C$. Because of the corollary, each geometric fiber of $\overline
f$ is isomorphic to ${\bf P}^1$, so $X$ is a ${\bf P}^1$-bundle on ${\bf
A}^1$. If we call $E$ the complement of ${\bf A}^2$ in $X$, with its reduced
scheme structure, then the projection from $E$ onto ${\bf A}^1$ is an
isomorphism. Hence there is an isomorphism $\phi$ of ${\bf P}^1\times {\bf
A}^1$ with $X$ carrying ${\bf A}^1\times {\infty}$ into $E$, and such that
$\overline f \circ \phi\colon {\bf P}^1\times {\bf A}^1\to {\bf A}^1$ is the
second projection. The restriction of $\phi$ to ${\bf A}^2$ carries the line
with equation $y = 0$ into $C$, and this proves the theorem.\endproof

This can be made to work in positive characteristic. Let
us fix an algebraically closed field $\kappa$, and call
$p$ be the characteristic exponent of $\kappa$, namely the characteristic of
$\kappa$ if this is positive, and 1 otherwise.

Consider a primitive polynomial in two variables, i.e., a morphism $f\colon
{\bf A}^2\to {\bf A}^1$ defined over $\kappa$ with integral general fiber;
as before, $f$ has a minimal normal compactification $\overline f \colon X
\to {\bf A}^1$. Define $E_1, \ldots, E_r$, $e_1, \ldots, e_r$ and $\delta_1,
\ldots, \delta_r$ as before. Then we can not conclude that the $\delta_ie_i$
are relatively prime; however, we have the following.

\proclaim Theorem 2.  Each of the $E_1,\ldots,E_r$ is isomorphic to
${\bf A}^1$, and they are pairwise disjoint.

Furthermore the greatest common divisor of the products
$\delta_1e_1,\ldots, \delta_re_r$ is a power of
$p$ and divides the degree of $f$.
\endproclaim

We still get the corollary, in the following form.

\proclaim Corollary 2.  Assume that the degree of $f$ is prime to the
$p$, and that $X$ has only one component at infinity.
Then all the fibers of $\overline f$ are integral and smooth at
infinity.

In particular, this happens when one of the fibers of $f$ has
only one branch at infinity.\endproclaim

Remarkably, using a different technique one can prove that $X$ is smooth when
it has only one component at infinity, without assuming that the degree of
$f$ is prime to $p$. Unfortunately, I do not have any interesting
application of this.

{}From Corollary~2 we get a new proof of the Suzuki--Abhyankar--Moh
theorem over any perfect field.

\proclaim The Suzuki--Abhyankar--Moh theorem over a perfect field.  Any
embedding of ${\bf A}^1$ into ${\bf A}^2$ defined over a perfect field, whose
degree is relatively prime to the characteristic of $\kappa$ is equivalent
to the standard embedding $t\mapsto(t,0)$.
\endproclaim

\proofof Theorem 1.  Recall that $\Gamma$ is the
closure of the graph of $f$ in ${\bf P}^2\times {\bf P}^1$, $X'$ its
normalization,
$f'$ and $\pi$ the projections of $X'$ onto ${\bf P}^1$ and ${\bf P}^2$,
respectively. Let $L = {\bf P}^2\setminus{\bf A}^2$ be the line at infinity,
$L'$ its proper transform in $X'$.

Let $E'_i$ be the closure of $E_i$ in $X'$: the first statement of the
theorem is a consequence of the following fact.

\proclaim Lemma 1.  The curves $L'$ and $E'_i$, for each $i = 1,\ldots,r$,
are isomorphic to ${\bf P}^1$, and any two of them do not intersect in more
than one point. Furthermore, if $E'_i$ and $E'_j$, with $i \neq j$, intersect
in a closed point $p \in X'$, then $p\in L'$.
\endproclaim

Assuming Lemma 1, and keeping in mind that that $L'$ is the fiber of $f'$
over the point at infinity $\infty\in {\bf P}^1(\kappa)$, we see that each
of the $E'_i$ can have only one point over $\infty$, and therefore the
inverse image $E_i$ of ${\bf A}^1$ in $E'_i$ is isomorphic to ${\bf A}^1$.
Also {}from Lemma 1 we get that the $E_i$ do not intersect.

\proof The natural morphism $\pi\colon X'\to {\bf P}^2$ is birational and
${\bf P}^2$ is smooth, so ${\rm R}^1\pi_*{\cal O}_{X'} = 0$.   Let
$\widetilde L = \pi^{-1}(L)_{\rm red}$. Since ${\cal O}_{\widetilde L}$ is a
quotient of ${\cal O}_{X'}$, so  ${\rm R}^1\pi_*{\cal O}_{\widetilde L} = 0$.
We have
$\pi_*{\cal O}_{\widetilde L} = {\cal O}_L$, so {}from the Leray spectral
sequence
$$
E_2^{ij} = \H^i({\bf P}^1,{\rm R}^j\pi_*{\cal O}_{\widetilde L})
\Longrightarrow \H^{i+j}(\widetilde L,{\cal O})
$$
we get that $\H^1(\widetilde L,{\cal O}) =
0$. If $Z$ is subscheme of $\widetilde L$, the sheaf ${\cal O}_Z$ is a
quotient of ${\cal O}_{\pi^{-1}(L)}$, and if ${\cal I}$ is the ideal of $Z$
in $\pi^{-1}(L)$ we have $\H^2(\pi^{-1}(L),{\cal I}) = 0$, hence
$\H^1(Z,{\cal O}) = 0$. This in particular applies to any of the curves $L'$
and $E'_i$. Any integral projective curve with arithmetic genus 0 is
isomorphic to ${\bf P}^1$.

Also, if $C_1$ and $C_2$ are two of these curve, {}from the fact
that $\H^1(C_1\cup C_2,{\cal O}) = 0$ we see that $C_1$ and $C_2$ have at
most one common point. Analogously, the fact that $\H^1(L'\cup E'_i\cup
E'_j,{\cal O}) = 0$ implies that $E'_i$ and $E'_j$ cannot meet outside of
$L'$, because
$L'$ meets both $E'_i$ and $E'_j$.\endproof

Now consider the group $\pic X$ of Cartier divisors on $X$, and the natural
map $\pic X \to \cl X$ into the group of Weil divisors. Since ${\bf A}^2$ is
factorial, and all of its invertible regular functions are constant, it
follows that $\cl X$ is a free abelian group with basis $E_1$, \dots,~$E_r$.
Since the map $\pic X \to \cl X$ is injective, because $X$ is normal, this
proves the following.

\proclaim Lemma 2. The group $\pic X$ is free, with basis
$\delta_1E_1$, \dots,~$\delta_rE_r$.
\endproclaim

The fact that the
$\delta_ie_i$ are relatively prime is easily proved, after having
established the following two facts.

\proclaim Lemma 3. The first Chern class map $\pic X \to \H^2(X, {\bf Z})$
is an isomorphism.
\endproclaim

\proclaim Lemma 4. Let $C$ be a general fiber of $\overline f$. The the
restriction map $\H^2(X, {\bf Z}) \to \H^2(C, {\bf Z})$ is
surjective. 
\endproclaim

In fact, the restriction of $\delta_iE_i$ to $C$ has degree $\delta_ie_i$;
the three lemmas together imply that the restriction of the $\delta_iE_i$
generate $\H^2(C, {\bf Z}) = {\bf Z}$, hence that 1 is a linear
combination of the $\delta_ie_i$.

There remains to give proofs of the last two lemmas; both are rather formal.

\proofof Lemma 3. Let $\rho\colon \widetilde X \to X$ be a resolution of the
singularities of $X$, $F_1, \ldots, F_s$ the exceptional divisors,
$\widetilde E_i$ the proper transforms of the $E_i$. Then the complement of
the $F_j$ and the $\widetilde E_i$ in $\widetilde X$ is ${\bf A}^2$;
therefore the Picard group of $\widetilde X$ is freely generated by the
$F_j$ and the $\widetilde E_i$. Likewise, $\H^2(\widetilde X, {\bf Z})$ is
freely generated by the cohomology classes of the $E_i$ and $F_j$; so the
first Chern class map $\pic \widetilde X \to \H^2( \widetilde X, {\bf Z})$
is an isomorphism.

The pullback map $\pic X \to \pic \widetilde X$ is clearly
injective, and a divisor class in $\pic \widetilde X$ is in the
image of $\pic X$ if and only if its restriction to each of the $F_j$ has
degree 0. The reason is that $X$ has rational singularities ([Li],
Theorem~12.1

Now take cohomology. We have that ${\rm R}^1\rho_* {\bf Z}_{ \widetilde X} =
0$, while $\rho_* {\bf Z}_{ \widetilde X} = {\bf Z}$, and ${\rm R}^2\rho_*
{\bf Z}_{\widetilde X}$ is a sheaf concentrated in the singular points of
$X$, whose stalk over $p \in X$ is a direct sum of one copy of {\bf
Z} for each exceptional divisor over $p$. By considering the Leray spectral
sequence of the map $\rho\colon\widetilde X \to X$, one deduces that the
restriction map $\H^2(X, {\bf Z}) \to \H^2( \widetilde X, {\bf Z})$ is
injective, and its image consists exactly of the classes in $\H^2(
\widetilde X, {\bf Z})$ which have degree 0 on each $F_j$.

By putting these two statements together, we see that $\pic X$ and $\H^2(X,
{\bf Z})$ are identified with two subgroups of $\pic \widetilde X$ and
$\H^2(\widetilde X, {\bf Z})$ which correspond under the isomorphism $\pic
\widetilde X \to \H^2(\widetilde X, {\bf Z})$ given by the first Chern
class. This proves Lemma~3.\endproof

\proofof Lemma~4. Consider the Leray spectral sequence
$$
E_2^{ij} = \H^i\bigl({\bf A}^1,{\rm R}^j\overline f_*{\bf Z}_X\bigr)
\Longrightarrow \H^{i+j}(X,{\bf Z});
$$
since $\H^2\bigl({\bf A}^1,{\rm R}^1\overline f_*{\bf Z}_X\bigr) =
0$, because ${\bf A}^1$ is an affine curve and ${\rm R}^1\overline f_*{\bf
Z}_X$ a constructible sheaf, we get that the map
$$
\H^2(X,{\bf Z})\to \H^0({\bf A}^1,{\rm R}^2\overline f_*{\bf Z}_X)
$$
is surjective. Now take the trace map
$$
\tr: {\rm R}^2\overline f_*{\bf Z}_X \to {\bf Z}_{{\bf A}^1}.
$$
Because the general fiber of $\overline f $ is integral, the trace
map is generically an isomorphism. Let $F = \overline f ^{-1}(t)$ be a fiber
of $\overline f $ over a closed point $t\in {\bf A}^1(\kappa)$,
$F_1,\ldots, F_s$ the irreducible components of $F$, $m_1,\ldots,m_s$ the
lengths of the local rings of $F$ at $F_1,\ldots,F_s$. By proper base
change the stalk $\bigl({\rm R}^2\overline f_*{\bf Z}_X)_t$ is canonically
isomorphic to
$$
\H^2(F,{\bf Z}) \simeq \bigoplus_{i=1}^s \H^2(F_i,{\bf Z}) \simeq 
{\bf Z}^s;
$$
with this identification, the trace map on the stalks
over $t\in {\bf A}^1(\kappa)$ is identified with the map {}from ${\bf Z}^s$
to {\bf Z} that sends $(k_1,\ldots,k_s)$ to
$k_1m_1+\cdots +k_sm_s$. But $m_1,\cdots,m_s$ are relatively prime, because
$f\colon S\to {\bf A}^1$ does not have multiple fibers, so the trace map is
surjective, and its kernel is concentrated on a finite number of points. By
taking global sections we see that the global trace map
$$
\tr: \H^2(X,{\bf Z})\to {\bf Z}
$$
is surjective. But $\tr: \H^2(X,{\bf Z})\to {\bf Z}$ coincides
with the restriction map $\H^2(X,{\bf Z})\to \H^2(C,{\bf Z})
\simeq  {\bf Z}$. Hence this restriction map is surjective. This
proves the lemma, and hence the theorem.\endproof

\proclaimr Note. {}From the spectral sequence of the map $X \to {\bf A}^1$
one deduces that $\H^3(X, {\bf Z}) = 0$; furthermore, {}from the spectral
sequence of a resolution $\widetilde X \to X$ one sees that the restriction
map $\H^2( \widetilde X, {\bf Z}) \to \H^2(F, {\bf Z})$ is surjective. {}From
this one can deduce that the class of $E_i$ generates
the product
$\prod_{p \in E_i}\cl
\widehat {\cal O}_{X,p}$; this means that $\delta_i$ can also be defined as
the product of the orders of the group $\cl\widehat {\cal O}_{X,p}$ for $p
\in E_i$.
\endproclaim

To prove Theorem~2 one follows the steps in the proof of Theorem~1,
subsituting \'etale cohomology with ${\bf Z}_\ell$ coefficients to classical
cohomology, where
$\ell$ is a prime different {}from the characteristic of
$\kappa$; in this way one shows that $\ell$ does not divide the greatest
common divisor of the $e_i \delta_i$. We leave the details to the interested
reader. The only thing that does not follow is that the greatest common
divisor of the
$\delta_ie_i$ divides the degree $d$ of $f$.

To show this, call $C$ the closure in ${\bf P}^2$
of a general fiber of $f$, $C'$ the proper transform of $C$ in $X'$. Then
$C'$ is a general fiber of $f'$, hence it is a Cartier divisor on $X'$; the
intersection number
$(C'\cdot L')$ is 0, and $(C'\cdot E'_i)= e_i$ for each $i = 1,\ldots,r$. We
have a decomposition of $\pi^*(L)$ as a Weil divisor
$$
\pi^*[L] = [L'] + \sum_{i=1}^r m_i E_i
$$
for certain positive integers
$m_1,\ldots,m_r$. Since the restriction of the $m_iE_i$ to $X$ must be a
Cartier divisor, we see that $\delta_i$ divides $m_i$, so we write
$$
\pi^*[L] = [L'] + \sum_{i=1}^r n_i \delta_iE_i.
$$
But
$$
d = (C\cdot L) = (C'\cdot \pi^*[L]) = (C'\cdot L') + \sum_{i=1}^r
n_i\delta_i (C'\cdot E'_i) = \sum_{i=1}^r n_i\delta_ie_i,
$$
by the projection formula, and this completes the proof.\endproof

\references

\w AM Abhyankar, S. S., Moh, T. T.:  Embeddings of the line in
the plane, J. Reine Angew. Math. {\bf 276}, 148--166 (1975).

\w AS Abhyankar, S. S., Singh, B.:  Embeddings of certain
curves, Am. J. of Math. {\bf 100}, 99--175 (1978).

\w AB Artal Bartolo, E.: Une d\'emonstration g\'eom\'etrique du th\'eor\`eme
d'Abhyankar--Moh, J. Reine Angew. Math. {\bf 464}, 97--108 (1995).

\w EN Eisenbud, D., Neumann, W. D.: Three-dimensional link theory and
invariants of plane curve singularities, Ann. Math. Studies {\bf101}, Princeton
University Press, Princeton N.Y. (1985).

\w Li Lipman, J.:  Rational singularities, with applications to
algebraic surfaces and unique factorization, Publications Math\'ematiques
I.H.E.S. {\bf 36},  195--280 (1969).

\w Se Segre, B.:  Forme differenziali e loro integrali, Docet, Roma
(1956).

\w Su Suzuki, M.: Proprietes topologiques des polynomes de deux
variables complexes et automorphismes alg\'ebrique de l'espace ${\bf C}^2$,
J. Math. Soc. Japan {\bf 26}, 241--257 (1974).

\bye